\documentclass[12pt]{article}
\usepackage{amssymb,amscd}
\usepackage{amsmath}
\usepackage{latexcad}

\newtheorem{conj}{Conjecture}[section]
\newtheorem{thm}[conj]{Theorem}

\newtheorem{cor}[conj]{Corollary}

\begin{document}
\title{\Large Equipartitions of a Mass in Boxes}
\author{{\Large Sini\v sa T. Vre\' cica}
\thanks{Supported by the
Serbian Ministry of Science, Grant 144026.}
\\ {University of Belgrade}
\\ {vrecica@matf.bg.ac.yu}} \maketitle \vskip 2cm

\begin{abstract}
The aim of this paper is to provide the sufficient condition for a
mass distribution in $\mathbb{R}^d$ to admit an equipartition with
a collection of hyperplanes some of which are parallel. The
results extend the previously obtained results for the
equipartitions with non-parallel hyperplanes. (See \cite{mvz} and
\cite{Ram}.)

The paper also serves as the illustration of the applicability and
the power of the methods of equivariant topology (more precisely,
equivariant index theory) in the problems of geometric
combinatorics.
\end{abstract}

\section{Introduction}

Any collection of $k$ hyperplanes in $\mathbb{R}^d$ determine a
partition of this Euclidean space (and any mass distribution in
it) into $2^k$ hyperorthants (defined as the intersections of the
appropriate half-spaces). Given a family of $j$ mass distributions
in $\mathbb{R}^d$, we say that a collection of $k$ hyperplanes
forms a equipartition of these $j$ mass distributions if each
hyperorthant contains exactly $\frac 1{2^k}$ of each of the given
mass distributions.

The question when every family of $j$ mass distributions in
$\mathbb{R}^d$ admit an equipartition by some collection of $k$
hyperplanes, is known as the equipartition problem, and it was
formulated by B. Gr\" unbaum in 1960 (see \cite{Gr}).

It attracted a lot of attention and some answers to this problem
are already obtained in \cite{Had}. Very thorough treatment of
this question is presented in \cite{Ram}, where more complete
results are obtained. However, the question remains unsettled in
general, and is still considered as an important and difficult
question in the area of discrete and computational
geometry/topology. The most recent and the most complete answers
to this question are given in \cite{mvz}. The most important open
question is whether every mass distribution in $\mathbb{R}^4$
admits an equipartition by $4$ hyperplanes in $16$ hyperorthants.

In this paper we treat the related problem of equipartition of a
mass distribution by a family of hyperplanes in a special
position, namely by a collection of parallel hyperplanes and one
or more additional non-parallel hyperplanes. Since such a
collection of hyperplanes divides $\mathbb{R}^d$ in the box-like
regions (or boxes), we will refer to this question as to the
question of equipartition of a mass distribution in boxes.

We obtain the general sufficient condition on the dimension $d$
and the number of parallel hyperplanes so that every mass
distribution in $\mathbb{R}^d$ admits such equipartition. (See
theorems \ref{2}, \ref{three}, \ref{threeo}, \ref{gen} and
\ref{geno}.) As the sample results illustrating the obtained
results we mention here the following three (see corollary
\ref{two}, corollary \ref{trio} and corollary \ref{tri}).
\bigskip

\noindent {\bf Claim 1.} Every mass distribution in the plane
admits an equipartition in $6$ boxes by two parallel lines and one
additional line not parallel to them.
\bigskip

\noindent {\bf Claim 2.} Any mass distribution in $\mathbb{R}^4$
could be equipartitioned in $12=3\times 2\times 2$ boxes by a
collection of $4$ hyperplanes two of which are parallel.
\bigskip

\noindent {\bf Claim 3.} Any mass distribution in $\mathbb{R}^8$
could be equipartitioned in $7\times 2\times 2$ boxes by a
collection of $6$ parallel hyperplanes and two additional
non-parallel hyperplanes.
\bigskip

Notice that the claim 2 is related to the (above mentioned) most
important open case of the original question of B. Gr\" unbaum. We
don't know whether for every mass distribution in $\mathbb{R}^4$,
there exists a $4$-tuple of hyperplanes equipartitioning this mass
distribution in $16$ hyperorthants. But, if we consider $4$-tuples
of hyperplanes two of which are parallel, they divide the space
and the mass distribution in $12$ boxes, and claim 2 shows that we
could always find such a $4$-tuple equipartitioning the measure.

For technical reasons, we treat separately the cases of even and
odd number of parallel hyperplanes. We first discuss the case of
even number of parallel hyperplanes, and then explain the
differences in the formulation and the proof of the odd case.

Throughout this paper, we work with the continuous mass
distributions with the positive measure of any open set in
$\mathbb{R}^d$. (A continuous mass distribution is a finite Borel
measure $\mu$ defined by the formula $\mu(A) = \int_A f \, d\mu$
for an integrable density function $f : {\mathbb R}^d\rightarrow
R$.) Because of that, the hyperplane orthogonal to some direction
and partitioning the given mass distribution in the given ratio is
unique. Using the limit argument, it is easy to extend the result
to all mass distributions which are weak limits of the mass
distributions satisfying the above properties. In particular, the
results are true for measurable sets and for finitely supported
measures.

We use the topological method in treating this question. More
precisely, we reduce the above question to the question of the
existence of an equivariant map. There is a number of ways to
treat the latter question, such as the use of characteristic
classes or the use of the obstruction theory. We find it most
convenient to use the index theory approach as formulated by E.
Fadell and S. Husseini in \cite{fh}.

The application of topological methods in combinatorics dates back
(at least) to 70's and the papers by L. Lov\' asz, I. B\' ar\' any
and others. The appearance of these ideas and their development
served as the starting point in the creation of the new subfield,
topological combinatorics.

\section{A short review of index theory}

For the reader's convenience, we present a very short review of
the ideal-valued cohomological index theory by E. Fadell and S.
Husseini. Given a finite group $G$, and a $G$-map $f : X
\rightarrow Y$ between $G$-spaces $X$ and $Y$, we could map these
spaces to the one-point space $\{\ast \}$ and obtain a commutative
diagram of $G$-spaces and $G$-maps. Multiplying by the total space
$EG$ of the universal $G$-bundle $EG\rightarrow BG$, we obtain new
commutative diagram of $G$-spaces and $G$-maps. (We consider the
diagonal $G$-action on the product spaces.) Passing to the spaces
of orbits, we obtain the following commutative diagram of
continuous maps:

\begin{figure}[htb]
\special{em:linewidth 0.4pt} \unitlength 1.00mm
\linethickness{0.4pt}
\begin{picture}(92,34)
\thinlines \drawcenteredtext{44.0}{26.0}{$X\times_GEG$}
\drawcenteredtext{102.0}{26.0}{$Y\times_GEG$}
\drawcenteredtext{74.0}{6.0}{$BG$}
\drawvector{56.0}{26.0}{34.0}{1}{0}
\drawvector{46.0}{24.0}{22.0}{4}{-3}
\drawvector{102.0}{24.0}{22.0}{-4}{-3}
\drawcenteredtext{54.0}{14.0}{$\tilde p_1$}
\drawcenteredtext{94.0}{14.0}{$\tilde p_2$}
\drawcenteredtext{72.0}{28.0}{$\tilde f$}
 \end{picture}
\caption{}
\end{figure}

\noindent which induces the following commutative diagram in
cohomology:

\begin{figure}[htb]
\special{em:linewidth 0.4pt} \unitlength 1.00mm
\linethickness{0.4pt}
\begin{picture}(98,38)
\thinlines \drawcenteredtext{46.0}{28.0}{$H_G^{\ast}(X)$}
\drawcenteredtext{106.0}{28.0}{$H_G^{\ast}(Y)$}
\drawcenteredtext{76.0}{6.0}{$H^{\ast}(BG)$}
\drawvector{94.0}{28.0}{36.0}{-1}{0}
\drawvector{68.0}{10.0}{20.0}{-4}{3}
\drawvector{84.0}{10.0}{20.0}{4}{3}
\drawcenteredtext{76.0}{32.0}{$\tilde f^{\ast}$}
\drawcenteredtext{52.0}{16.0}{$\tilde p_1^{\ast}$}
\drawcenteredtext{100.0}{16.0}{$\tilde p_2^{\ast}$}
\end{picture}
\caption{}
\end{figure}

The kernels of the maps $\tilde p_1^{\ast}$ and $\tilde
p_2^{\ast}$ are the ideals in the cohomology ring of the
classifying space of the group $G$, and they are called indices
and denoted by $\mbox{Ind}_G X$ and $\mbox{Ind}_G Y$ respectively.

The commutativity of the above diagram implies the relation
$\mbox{Ind}_G Y \subseteq \mbox{Ind}_G X$. If we could prove that
this inclusion relation is not satisfied, we would obtain a
contradiction proving that a $G$-equivariant map $f : X\rightarrow
Y$ does not exist.

We refer the reader to the original paper \cite{fh} for additional
properties of the index and some basic computation. Some other
computations, needed in this paper could be found in \cite{z},
Corollary 2.12 and Proposition 2.7.

\section{The problem and the results - 2 directions}

In order to describe the method and develop the intuition in a
more acceptable way, we choose to treat the simplest particular
case first.

So, in this section we treat the question of equipartition of a
mass distribution in $\mathbb{R}^d$ by a collection of parallel
hyperplanes and by one additional hyperplane (not parallel to
them). We show that the greatest number of parallel hyperplanes
for which such an equipartition always exists (for every mass
distribution in $\mathbb{R}^d$) is $2d-2$. The same result will be
true also for $2d-3$ parallel hyperplanes in $\mathbb{R}^d$, but
the equipartition with $2d-1$ parallel hyperplanes are always
possible only in $\mathbb{R}^{d+1}$. The obtained result is the
best possible in the sense that for greater number of hyperplanes,
the corresponding equivariant mapping exists.

\begin{thm}
\label{2} For every mass distribution in $\mathbb{R}^d$ there is a
collection of $2d-2$ parallel hyperplanes and one additional
hyperplane dividing $\mathbb{R}^d$ in $4d-2$ boxes containing the
same amount of the mass distribution.
\end{thm}

Especially, when $d=2$, we obtain the proof of the following
corollary.

\begin{cor}
\label{two} Every mass distribution in the plane admits an
equipartition in $6$ boxes by two parallel lines and one
additional line not parallel to them.\hfill $\blacksquare$
\end{cor}

The proof of the theorem \ref{2} (with the complete description of
the approach, needed also in the proofs of other results from this
paper) is contained in the following two subsections.

\subsection{Reduction} \label{red}

In this subsection we reduce the statement of the above theorem to
the topological statement.

For any mass distribution in $\mathbb{R}^d$ and any pair of
vectors $(u,v)\in S^{d-1}\times S^{d-1}$, let $H_1^u,
H_2^u,...,H_{2d-2}^u$ be the oriented hyperplanes orthogonal to
$u$, ordered in the direction of the vector $u$, and dividing
$\mathbb{R}^d$ into $2d-1$ regions each containing the same amount
(i.e. $\frac 1{2d-1}$) of the considered mass distribution. Also,
let $H^v$ be the oriented hyperplane orthogonal to $v$ dissecting
a mass distribution into two halfspaces containing the same amount
of the mass distribution.

These hyperplanes form $2(2d-1)$ boxes and the measure of these
boxes form a $2\times (2d-1)$ matrix of the form
$$\left(
\begin{array}{rrrr}
\alpha_0+\alpha_1 & \alpha_0+\alpha_2 & \dots &
\alpha_0+\alpha_{2d-1}\\
\alpha_0-\alpha_1 & \alpha_0-\alpha_2 & \dots &
\alpha_0-\alpha_{2d-1}
\end{array}  \right) $$

\noindent where $\alpha_0=\frac 1{2(2d-1)}$ and
$\alpha_1+\alpha_2+\cdots +\alpha_{2d-1}=0.$

So, we could identify the configuration space of our problem to be
the product of two spheres $S^{d-1}\times S^{d-1}$ and the test
space as the space of all $2\times (2d-1)$ matrices of the above
form. The group $\mathbb{Z}/2\oplus \mathbb{Z}/2$ acts naturally
both on configuration space and the test space (by the obvious
permutations). The test space could also be seen as the
$(2d-2)$-dimensional linear representation of the group
$\mathbb{Z}/2\oplus \mathbb{Z}/2$, which we denote by $V$. The
test map $f : S^{d-1}\times S^{d-1} \rightarrow V$, which maps
each pair of unit vectors to the measures of the corresponding
boxes, is easily seen to be $(\mathbb{Z}/2\oplus
\mathbb{Z}/2)$-equivariant.

Now, our problem is reduced to the topological claim that the
matrix
$$\left(
\begin{array}{rrrr}
\alpha_0 & \alpha_0 & \dots & \alpha_0\\
\alpha_0 & \alpha_0 & \dots & \alpha_0
\end{array}  \right) $$

\noindent (obtained when $\alpha_1=\cdots =\alpha_{2d-1}=0$)
belongs to the image of the test map $f$. Suppose, to the
contrary, this not to be the case. Then we would have a
$(\mathbb{Z}/2\oplus \mathbb{Z}/2)$-equivariant map $f :
S^{d-1}\times S^{d-1} \rightarrow S(V)$, where $S(V)$ denotes the
unit sphere of the representation space $V$. Finally, we reach a
contradiction (proving in this way our claim), by showing that
such equivariant map with the actions of our group described
formerly could not exist. In proving this we use the ideal valued
cohomological index theory of Fadell and Husseini.

\subsection{Computation}

We will use the approach described above to show that there is no
$(\mathbb{Z}/2\oplus \mathbb{Z}/2)$-equivariant map $f :
S^{d-1}\times S^{d-1} \rightarrow S(V)$, where $S(V)$ denotes the
unit sphere of the representation space $V$ described in the
section \ref{red}. So, we work with the group
$G=\mathbb{Z}/2\oplus \mathbb{Z}/2$, and it is well known that
$BG=B\mathbb{Z}/2\times B\mathbb{Z}/2=\mathbb{R}P^{\infty}\times
\mathbb{R}P^{\infty}$, and $H^{\ast}(BG;\mathbb{Z}/2)\cong
\mathbb{Z}/2[x,y]$, where $x$ and $y$ are the free generators of
this polynomial ring in dimension $1$ both.

The generators of the group $G=\mathbb{Z}/2\oplus \mathbb{Z}/2$
act by the antipodal action on the corresponding spheres in
$S^{d-1}\times S^{d-1}$. It is well known that in this case we
have $\mbox{Ind}_G(S^{d-1}\times S^{d-1})=(x^d,y^d)$, i.e. the
index is the ideal generated by the monomials $x^d$ and $y^d$.

Now we determine the index of the unit sphere in the
representation space $V$. We refer the reader to the survey
article \cite{z} (Corollary 2.12), for the necessary background
for the following computation. As we noticed in the section
\ref{red}, $V$ is the $(2d-2)$-dimensional representation which
could be described as the space of all $2\times (2d-1)$ matrices
of the form
$$\left(
\begin{array}{rrrr}
\alpha_0+\alpha_1 & \alpha_0+\alpha_2 & \dots &
\alpha_0+\alpha_{2d-1}\\
\alpha_0-\alpha_1 & \alpha_0-\alpha_2 & \dots &
\alpha_0-\alpha_{2d-1}
\end{array}  \right) $$

\noindent where $\alpha_0=\frac 1{2(2d-1)}$ and
$\alpha_1+\alpha_2+\cdots +\alpha_{2d-1}=0.$ The generator of the
first copy of $\mathbb{Z}/2$ acts on such matrices by permuting
the columns in the reverse order, i.e. by sending
$(\alpha_1,\alpha_2,...,\alpha_{2d-1})$ to
$(\alpha_{2d-1},\alpha_{2d-2},...,\alpha_1)$. The generator of the
second copy of $\mathbb{Z}/2$ acts by permuting two rows, i.e. by
sending each $\alpha_i$ to $-\alpha_i$. By the relation
$\alpha_1+\alpha_2+\cdots +\alpha_{2d-1}=0$, the element
$\alpha_d$ is determined by the remaining elements. To shorten the
notation we will subtract $\alpha_0$ from the entries of the above
mentioned matrix, and present the matrix in the form
$$\left(
\begin{array}{rrrr}
\alpha_1 & \alpha_2 & \dots & \alpha_{2d-1}\\
-\alpha_1 & -\alpha_2 & \dots & -\alpha_{2d-1}
\end{array}  \right) $$

So, $V$ is a $(2d-2)$-dimensional representation of the group
$\mathbb{Z}/2\oplus \mathbb{Z}/2$. The representation space $V$
splits in the sum of $2d-2$ invariant $1$-dimensional
representations. We present $d-1$ pairs of this $G$-invariant
$1$-dimensional representations. The $i$-th pair of this
representations form the matrices of the form:
$$\left(
\begin{array}{rrrrrrr}
\dots & \alpha & \dots & -2\alpha & \dots & \alpha & \dots \\
\dots & -\alpha & \dots & 2\alpha & \dots & -\alpha & \dots
\end{array}  \right) $$

\noindent and the matrices of the form:
$$\left(
\begin{array}{rrrrrrr}
\dots & \alpha & \dots & 0 & \dots & -\alpha & \dots \\
\dots & -\alpha & \dots & 0 & \dots & \alpha & \dots
\end{array}  \right) $$

Here, we write only the entries in the $i$-th, $d$-th, and
$(2d-i)$-th column of the matrix, while all other entries are $0$.

The generator of the first copy of $\mathbb{Z}/2$ acts on such
matrices by permuting the columns in the reverse order, and so it
acts trivially on the first mentioned $1$-dimensional subspace of
matrices and antipodally on the second. The generator of the
second copy of $\mathbb{Z}/2$ acts by permuting two rows, and so
it acts antipodally on both $1$ dimensional subspaces of matrices.

So, the index $\mbox{Ind}_G S(V)$ is the ideal in the polynomial
ring $\mathbb{Z}/2[x,y]$ generated by the polynomial
$(y(x+y))^{d-1}$. (Consult \cite{z}.)

All the summands of the polynomial
$$y^{d-1}(x+y)^{d-1}=
x^{d-1}y^{d-1}+{d-1\choose 1}x^{d-2}y^d+\cdots +y^{2d-2},$$

\noindent belong to the ideal $(x^d,y^d)$, except for the first
one $x^{d-1}y^{d-1}$. So, this polynomial is not contained in the
ideal $(x^d,y^d)$. This means that $\mbox{Ind}_G S(V)\nsubseteq
\mbox{Ind}_G (S^{d-1}\times S^{d-1})$. The considerations from the
previous subsection show that there is no $(\mathbb{Z}/2\oplus
\mathbb{Z}/2)$-equivariant mapping from $S^{d-1}\times S^{d-1}$ to
$S(V)$, which implies that every equivariant map from
$S^{d-1}\times S^{d-1}$ to the representation space $V$ maps some
pair of unit vectors $(u,v)$ to the matrix
$$\left(
\begin{array}{rrrr}
\alpha_0 & \alpha_0 & \dots & \alpha_0\\
\alpha_0 & \alpha_0 & \dots & \alpha_0
\end{array}  \right) $$

\noindent This completes the argument and proves our theorem.
\hfill $\blacksquare$
\bigskip

If we consider the case of odd number of parallel hyperplanes
($2d-1$ of them), the similar considerations would prove the
following theorem.

\begin{thm}
For every mass distribution in $\mathbb{R}^{d+1}$ there is a
collection of $2d-1$ parallel hyperplanes and one additional
hyperplane dividing it in $4d$ boxes containing the same amount of
the mass distribution. \hfill $\blacksquare$
\end{thm}

The major difference is that we are now faced with two central
columns, since the matrix has even number of columns ($2d$ of
them), and in very similar way we get that the index of the sphere
in the representation space is the ideal generated by the
polynomial $y^{d-1}(x+y)^d$. This polynomial belongs to the ideal
generated by monomials $x^d$ and $y^d$, but does not belong to the
ideal generated by monomials $x^{d+1}$ and $y^{d+1}$. Notice that
in the case $m=2$ the stronger result (in some sense) is obtained
for even number of hyperplanes. Namely, any mass distribution in
$\mathbb{R}^{d+1}$ could be equipartitioned also in $4d+2$ boxes
by some $2d$ parallel hyperplanes and one additional non-parallel
to them.

\section{The case of 3 directions}

In this section we generalize the result from the previous section
to the case of equipartition of a mass distribution in some
Euclidean space $\mathbb{R}^d$ by a collection of parallel
hyperplanes orthogonal to the direction $u$, and by two additional
hyperplanes not parallel neither to the first mentioned collection
nor one to each other. First we consider even number $2k$ of
parallel hyperplanes. In this case we treat the equipartition of a
mass distribution in $(2k+1)\times 2\times 2$ boxes. Our aim is to
determine the sufficient condition on the dimension $d$ and the
number $k$ so that every mass distribution in $\mathbb{R}^d$
admits such an equipartition.

Since we use the index theory again, we provide an algorithm to
decide the above question for a pair of numbers $d$ and $k$, which
reduces the question to the question whether some polynomial
belongs to some ideal in the polynomial algebra
$\mathbb{Z}/2[x_1,x_2,x_3]$ over $3$ variables.

Let us denote with $\mathbb{P}_3(x_1,x_2,x_3) =x_1x_2x_3 (x_1+x_2)
(x_1+x_3)(x_2+x_3)(x_1+x_2+x_3)$ the Dickson polynomial in $3$
variables. Over $\mathbb{Z}/2$, this Dickson polynomial has also
another description $\mathbb{P}_3(x_1,x_2,x_3)=\sum_{\sigma \in
S_3} x_{\sigma (1)}^4x_{\sigma (2)}^2x_{\sigma (3)}$, and as a
consequence we get:

\begin{thm}
\label{three}
Let
\[
{\mathbb P}_3 = {\rm Det}\left[
\begin{array}{ccc}
x_1 & x_1^2 & x_1^4\\

x_2 & x_2^2 & x_2^4\\

x_3 & x_3^2 & x_3^4
\end{array}\right] \in\mathbb{Z}/2[x_1,x_2,x_3]
\]
be a Dickson polynomial. Then every measure in $\mathbb{R}^d$
admits an equipartition by a collection of $2k$ parallel
hyperplanes and two additional non-parallel hyperplanes in
$(2k+1)\times 2\times 2$ boxes if
\[
(x_2+x_3)\left(\frac 1{x_1}{\mathbb P}_3\right)^k\notin
(x_1^d,x_2^d,x_3^d).
\]
\end{thm}

\medskip\noindent
{\bf Proof:} Again, for any mass distribution in $\mathbb{R}^d$
and any triple of vectors $(u,v,w)\in S^{d-1}\times S^{d-1}\times
S^{d-1}$, let $H_1^u, H_2^u,...,H_{2k}^u$ be the oriented
hyperplanes orthogonal to $u$, ordered in the direction of the
vector $u$, and dividing $\mathbb{R}^d$ into $2k+1$ regions each
containing the same amount (i.e. $\frac 1{2k+1}$) of the
considered mass distribution. Also, let $H^v$ and $H^w$ be the
oriented hyperplanes orthogonal to $v$ and $w$ respectively, each
dissecting a mass distribution into two halfspaces containing the
same amount of the mass distribution.

These hyperplanes form $(2k+1)\times 2\times 2$ boxes and the
measure of these boxes form a $3$-dimensional $(2k+1)\times
2\times 2$ matrix. We describe this matrix by its $2k+1$
two-dimensional "slices" which are $2\times 2$ matrices, and are
of the form
$$\left(
\begin{array}{rr}
\varrho +\alpha_i & \varrho +\beta_i\\
\varrho +\gamma_i & \varrho +\delta_i
\end{array}  \right) $$

\noindent ($i=1,2,...,2k+1$), where $\varrho =\frac 1{4(2k+1)}$,
and $\alpha_i +\beta_i +\gamma_i +\delta_i=0$ for every
$i=1,2,...,2k+1$. Also, the entries of this $3$-dimensional matrix
satisfy two additional relations (coming from the properties of
the hyperplanes $H^v$ and $H^w$), and those are $\sum_i (\alpha_i
+\beta_i)=0$ and $\sum_i (\alpha_i +\gamma_i)=0$.

In this case the configuration space of our problem is the product
of three spheres $S^{d-1}\times S^{d-1}\times S^{d-1}$ and the
test space is the space of all $(2k+1)\times 2\times 2$ matrices
of the above form. The group $\mathbb{Z}/2\oplus
\mathbb{Z}/2\oplus \mathbb{Z}/2$ acts naturally both on
configuration space and the test space. So, equivalently the test
space could be represented as a $(6k+1)$-dimensional linear
representation $V$ of the group $\mathbb{Z}/2\oplus
\mathbb{Z}/2\oplus \mathbb{Z}/2$.

The test map $f : S^{d-1}\times S^{d-1}\times S^{d-1} \rightarrow
V$, which maps each triple of unit vectors to the measures of the
corresponding boxes, is easily seen to be $(\mathbb{Z}/2\oplus
\mathbb{Z}/2\oplus \mathbb{Z}/2)$-equivariant.

Now, our problem is reduced to the topological claim that the
matrix with all entries equal to $\varrho=\frac 1{4(2k+1)}$
(obtained when $\alpha_i=\beta_i=\gamma_i=\delta_i=0$ for every
$i=1,2,...,2k+1$) belongs to the image of the test map $f$.
Suppose, to the contrary, this not to be the case. Then we would
have a $(\mathbb{Z}/2\oplus \mathbb{Z}/2\oplus
\mathbb{Z}/2)$-equivariant map $f : S^{d-1}\times S^{d-1}\times
S^{d-1}\rightarrow S(V)$, where $S(V)$ denotes the unit sphere of
the representation space $V$.

As in the previous case, we reach a contradiction (proving in this
way our claim), by showing that such equivariant map with the
described actions of our group could not exist.

In this case we have the group $G=\mathbb{Z}/2\oplus
\mathbb{Z}/2\oplus \mathbb{Z}/2$, and it is well known that
$H^{\ast}(BG;\mathbb{Z}/2)\cong \mathbb{Z}/2[x_1,x_2,x_3]$, where
$x_1$, $x_2$, and $x_3$ are the free generators of this polynomial
ring in dimension $1$ all.

The generators of the group $G$ act by the antipodal action on the
corresponding spheres in $S^{d-1}\times S^{d-1}\times S^{d-1}$. It
is well known that in this case we have
$\mbox{Ind}_G(S^{d-1}\times S^{d-1}\times
S^{d-1})=(x_1^d,x_2^d,x_3^d)$, i.e. the index is the ideal
generated by the monomials $x_1^d$, $x_2^d$, and $x_3^d$.

Now we determine the index of the unit sphere in the
representation space $V$. The generator of the first copy of
$\mathbb{Z}/2$ permutes the "slices" of the matrix in the reverse
order, i.e. by sending $\alpha_i,\beta_i,\gamma_i,\delta_i$ to
$\alpha_{2k+2-i},\beta_{2k+2-i},\gamma_{2k+2-i},\delta_{2k+2-i}$.
The generator of the second copy of $\mathbb{Z}/2$ acts by
permuting two $2$-dimensional "rows" of the matrix, i.e. by
sending each $\alpha_i$ and $\beta_i$ to $\gamma_i$ and $\delta_i$
respectively. The generator of the third copy of $\mathbb{Z}/2$
also acts by permuting two $2$-dimensional "rows" of the matrix,
i.e. by sending each $\alpha_i$ and $\gamma_i$ to $\beta_i$ and
$\delta_i$ respectively.

To shorten the notation we will subtract $\varrho$ from the
entries of the above mentioned matrix, and present the "slices" of
the matrix in the form
$$\left(
\begin{array}{rr}
\alpha_i & \beta_i\\
\gamma_i & \delta_i
\end{array}  \right) $$

The representation space $V$ splits in the sum of $6k+1$
$G$-invariant $1$-dimensional representations. We present here $k$
$6$-tuples of this $G$-invariant $1$-dimensional representations
and additionally the last one. The $i$-th $6$-tuple of this
representations form the matrices whose $i$-th, $(k+1)$-th, and
$(2k+2-i)$-th "slices" are of the following forms (the "slices"
will be separated by the vertical lines, remember that the entries
in the remaining "slices" are all $0$):
$$\left(
\begin{array}{r|r|r|r|r|r|r}
\dots &
\begin{array}{rr}
\lambda & \lambda\\
-\lambda & -\lambda
\end{array} &
\dots &
\begin{array}{rr}
-2\lambda & -2\lambda\\
2\lambda & 2\lambda
\end{array} &
\dots &
\begin{array}{rr}
\lambda & \lambda\\
-\lambda & -\lambda
\end{array} &
\dots
\end{array}  \right)$$

$$\left(
\begin{array}{r|r|r|r|r|r|r}
\dots &
\begin{array}{rr}
\lambda & \lambda\\
-\lambda & -\lambda
\end{array} &
\dots &
\begin{array}{rr}
0 & 0\\
0 & 0
\end{array} &
\dots &
\begin{array}{rr}
-\lambda & -\lambda\\
\lambda & \lambda
\end{array} &
\dots
\end{array}  \right)$$

$$\left(
\begin{array}{r|r|r|r|r|r|r}
\dots &
\begin{array}{rr}
\lambda & -\lambda\\
\lambda & -\lambda
\end{array} &
\dots &
\begin{array}{rr}
-2\lambda & 2\lambda\\
-2\lambda & 2\lambda
\end{array} &
\dots &
\begin{array}{rr}
\lambda & -\lambda\\
\lambda & -\lambda
\end{array} &
\dots
\end{array}  \right)$$

$$\left(
\begin{array}{r|r|r|r|r|r|r}
\dots &
\begin{array}{rr}
\lambda & -\lambda\\
\lambda & -\lambda
\end{array} &
\dots &
\begin{array}{rr}
0 & 0\\
0 & 0
\end{array} &
\dots &
\begin{array}{rr}
-\lambda & \lambda\\
-\lambda & \lambda
\end{array} &
\dots
\end{array}  \right)$$

$$\left(
\begin{array}{r|r|r|r|r|r|r}
\dots &
\begin{array}{rr}
\lambda & -\lambda\\
-\lambda & \lambda
\end{array} &
\dots &
\begin{array}{rr}
0 & 0\\
0 & 0
\end{array} &
\dots &
\begin{array}{rr}
\lambda & -\lambda\\
-\lambda & \lambda
\end{array} &
\dots
\end{array}  \right)$$

$$\left(
\begin{array}{r|r|r|r|r|r|r}
\dots &
\begin{array}{rr}
\lambda & -\lambda\\
-\lambda & \lambda
\end{array} &
\dots &
\begin{array}{rr}
0 & 0\\
0 & 0
\end{array} &
\dots &
\begin{array}{rr}
-\lambda & \lambda\\
\lambda & -\lambda
\end{array} &
\dots
\end{array}  \right)
$$

The additional invariant $1$-dimensional representation have
non-zero entries only in the $(k+1)$-th "slice" and is of the form
$$\left(
\begin{array}{r|r|r}
\dots &
\begin{array}{rr}
\lambda & -\lambda\\
-\lambda & \lambda
\end{array} &
\dots \end{array} \right)
$$

According to the described action of the generators of the group
$G$ on these $3$-dimensional matrices, we see that the indices of
the $1$-dimensional representations from the $6$-tuples are
generated by the polynomials $x_2$, $x_1+x_2$, $x_3$, $x_1+x_3$,
$x_2+x_3$, $x_1+x_2+x_3$ (in this order), and the index of the
additional representation is $x_2+x_3$. So, by \cite{z} (Corollary
2.12), the index of the test space $\mbox{Ind}_G S(V)$ is the
ideal in the polynomial ring $\mathbb{Z}/2[x_1,x_2,x_3]$ generated
by the polynomial
$$\left(x_2(x_1+x_2)x_3(x_1+x_3)(x_2+x_3)(x_1+x_2+x_3)\right)^k(x_2+x_3)=
(x_2+x_3)\left(\frac 1{x_1}{\mathbb P}_3\right)^k$$

\noindent and the result follows. \hfill $\blacksquare$
\bigskip

It is easy now to determine the smallest dimension $d$ for any
$k$, for which our method provides the answer to the question. By
the properties of the binomial coefficients over $\mathbb{Z}/2$,
the best estimate is obtained when $k$ is a little bit smaller
than some power of $2$. For example, for $k=1$ we get $d=4$, and
for $k=3$ we get $d=8$.

\begin{cor}
\label{trio} Any mass distribution in $\mathbb{R}^4$ could be
equipartitioned in $12=3\times 2\times 2$ boxes by a collection of
$4$ hyperplanes two of which are parallel. \hfill $\blacksquare$
\end{cor}

\begin{cor}
\label{tri}
Any mass distribution in $\mathbb{R}^8$ could be
equipartitioned in $7\times 2\times 2$ boxes by a collection of
$6$ parallel hyperplanes and two additional non-parallel
hyperplanes.
\end{cor}

\medskip\noindent
{\bf Proof:} We will show that the coefficients in the polynomial
$(x_2+x_3)\left(\frac 1{x_1}{\mathbb P}_3\right)^3$ multiplying
the monomials $x_1^7x_2^7x_3^5$ and $x_1^7x_2^5x_3^7$ are
non-trivial. Since these monomials do not belong to the ideal
generated by the monomials $x_1^8$, $x_2^8$ and $x_3^8$, the
corollary follows.

The third power of the sum of some monomials (over $\mathbb{Z}/2$)
has non-zero coefficient multiplying the third power of the
monomials and the product of the square of some monomial with some
other monomial. It is easy to verify that there is only one way to
get $x_1^7x_2^7x_3^5$ (e.g.) in the expression
$(x_2+x_3)\left(\frac 1{x_1}{\mathbb P}_3\right)^3$ and that is
$x_2\cdot \left(\frac 1{x_1}\right)^3\cdot
\left(x_1^4x_3^2x_2\right)^2\cdot \left(x_2^4x_1^2x_3\right)$. So
the coefficient multiplying $x_1^7x_2^7x_3^5$ is non-zero, and we
are done. \hfill $\blacksquare$
\bigskip

We believe that in these cases (as in the case $m=2$), the
obtained results are the best possible. At least, for the greater
values of dimension $d$, the considered equivariant mapping
exists. Namely, if we consider $l$ parallel hyperplanes and two
additional non-parallel to them, the corresponding equivariant
mapping maps $S^{d-1}\times S^{d-1}\times S^{d-1}$ to the unit
sphere in the representation space which is of dimension $3l$. So,
in order for this equivariant mapping not to exist, we have to
have $l\leq d-2$.

Let us now turn to the case od odd number od parallel hyperplanes
and two additional hyperplanes. In almost the same way we obtain:

\begin{thm}
\label{threeo} Let
\[
{\mathbb P}_3 = {\rm Det}\left[
\begin{array}{ccc}
x_1 & x_1^2 & x_1^4\\

x_2 & x_2^2 & x_2^4\\

x_3 & x_3^2 & x_3^4
\end{array}\right] \in\mathbb{Z}/2[x_1,x_2,x_3]
\]
be a Dickson polynomial. Then every measure in $\mathbb{R}^d$
admits an equipartition by a collection of $2k+1$ parallel
hyperplanes and two additional non-parallel hyperplanes in
$(2k+2)\times 2\times 2$ boxes if
\[
\frac 1{x_2x_3}\left(\frac 1{x_1}{\mathbb P}_3\right)^{k+1}\notin
(x_1^d,x_2^d,x_3^d).
\]
\end{thm}

\medskip\noindent
{\bf Proof:} Let us only explain the differences. We have the
$3$-dimensional matrix with even number of slices and we have the
central {\bf pair} of slices. This central pair forms a $2\times
2\times 2$ matrix subject to $4$ relations. The corresponding
representation space splits in the sum of $4$ one-dimensional
representations, those from the $6$-tuple in the proof of theorem
\ref{three} having only zeros in the central slice. Therefore, the
corresponding index of this representation space is:

$$\left(\frac 1{x_1}{\mathbb P}_3\right)^k \cdot
(x_1+x_2)(x_1+x_3)(x_2+x_3)(x_1+x_2+x_3) = \frac
1{x_2x_3}\left(\frac 1{x_1}{\mathbb P}_3\right)^{k+1}$$

\noindent The result follows. \hfill $\blacksquare$
\bigskip

An easy algebraic calculation provides us with the following table
in which we describe, for small numbers $l$ (being even or odd) of
parallel hyperplanes, the smallest dimension $d$ of the Euclidean
space in which the equipartition with that many parallel
hyperplanes and two additional non-parallel to them, is always
possible.
\medskip

\begin{center}
\begin{tabular}{|c||c|c|c|c|c|c|c|} \hline
$l$ & 2 & 3 - 4 & 5 - 6 & 7 - 10 & 11 - 12 & 13 - 14 & 15 - 22 \\
\hline

$d$ & 4 & 7 & 8 & 13 & 15 & 16 & 25
\\ \hline
\end{tabular}
\end{center}
\medskip

Reading this table in the other direction, we see that in
$\mathbb{R}^4$ (and in $\mathbb{R}^5$ and $\mathbb{R}^6$) the
equipartition is always possible with $2$ parallel hyperplanes
(and two non-parallel, which we do not mention any further), in
$\mathbb{R}^7$ with $4$, in $\mathbb{R}^8$ (and up to
$\mathbb{R}^{12}$)  with $6$, in $\mathbb{R}^{13}$ (and in
$\mathbb{R}^{14}$) with $10$, in $\mathbb{R}^{15}$ with $12$, in
$\mathbb{R}^{16}$ (and up to $\mathbb{R}^{24}$) with $14$, in
$\mathbb{R}^{25}$ with $22$, and so on.

Notice that again for $l=14$ we get $d=16$ which is the best
possible in the same sense as above. Notice also that in all these
examples, due to the arithmetic reasons, the resulting Euclidean
space has the same dimension for the odd number $2k-1$ of parallel
hyperplanes and for the next even number $2k$ of parallel
hyperplanes.

\section{The general case}

It is obvious how to generalize these statements to the case of
more than $3$ directions, to obtain the complete algorithm for the
determination of the dimension $d$ so that any mass distribution
in $\mathbb{R}^d$ admits an equipartition in boxes. Without going
into details, we formulate the results obtained for the case of
$m$ directions. The reader could modify the above argument to
provide the proof for this statement.

Again, we formulate two separate statements, one for the case of
even and the other for the case of odd number of parallel
hyperplanes. Similarly to the previous case, with ${\mathbb
P}_m(x_1,...,x_m)$ we denote the Dickson polynomial in $m$
variables. Again, it is the product of all linear combinations of
these variables. Over $\mathbb{Z}/2$ it could also be described by
${\mathbb P}_m(x_1,...,x_m)= \sum_{\sigma \in S_m} x_{\sigma
(1)}^{2^{m-1}}\cdots x_{\sigma (m)}$. The Dickson polynomial
${\mathbb P}_{m-1}$ mentioned below will be in $m-1$ variables
$x_2,...,x_m$.

\begin{thm}
\label{gen}
Let
\[
{\mathbb P}_m = {\rm Det}\left[
\begin{array}{cccc}
x_1 & x_1^2 & \dots & x_1^{2^{m-1}}\\

x_2 & x_2^2 & \dots & x_2^{2^{m-1}}\\

\vdots & \vdots & \ddots & \vdots\\

x_m & x_m^2 & \dots & x_m^{2^{m-1}}
\end{array}\right] \in\mathbb{Z}/2[x_1,x_2,...,x_m]
\]
be a Dickson polynomial. Then every measure in $\mathbb{R}^d$
admits an equipartition by a collection of $2k$ parallel
hyperplanes and $m-1$ additional non-parallel hyperplanes in
$(2k+1)\times 2\times \cdots \times 2$ boxes if
\[
\frac 1{x_2x_3\cdots x_m}{\mathbb P}_{m-1}(x_2,...,x_m)\left(\frac
1{x_1}{\mathbb P}_m\right)^k\notin (x_1^d,x_2^d,...,x_m^d).
\]
\hfill $\blacksquare$
\end{thm}
\bigskip

In the case of odd number of parallel hyperplanes, we have the
following statement.

\begin{thm}
\label{geno} Let
\[
{\mathbb P}_m = {\rm Det}\left[
\begin{array}{cccc}
x_1 & x_1^2 & \dots & x_1^{2^{m-1}}\\

x_2 & x_2^2 & \dots & x_2^{2^{m-1}}\\

\vdots & \vdots & \ddots & \vdots\\

x_m & x_m^2 & \dots & x_m^{2^{m-1}}
\end{array}\right] \in\mathbb{Z}/2[x_1,x_2,...,x_m]
\]
be a Dickson polynomial. Then every measure in $\mathbb{R}^d$
admits an equipartition by a collection of $2k+1$ parallel
hyperplanes and $m-1$ additional non-parallel hyperplanes in
$(2k+2)\times 2\times \cdots \times 2$ boxes if
\[
\frac 1{x_2x_3\cdots x_m}\left(\frac 1{x_1}{\mathbb
P}_m\right)^{k+1}\notin (x_1^d,x_2^d,...,x_m^d).
\]
\hfill $\blacksquare$
\end{thm}

\section{Concluding remarks}

\subsection{Limitations of the method}

Our method does not provide the answer to the case when we
consider the collections of parallel hyperplanes in $2$ or more
directions. Namely, there are infinitely many fixed points of the
action of the group $G$ on the test space in these cases, and so
the equivariant map exists. The same is true if we consider the
case of more than one mass distribution.

\subsection{Acknowledgement}

This work initiated during the visit to the Mathematical Sciences
Research Institute, in Berkeley. The author wishes to thank MSRI
for their support and warm hospitality.

\vskip 1cm

\parbox{6cm}{Sini\v sa T. Vre\' cica \par Faculty of Mathematics
\par University of Belgrade
\par Studentski trg 16, P.O.B. 550 \par
\par 11000 Belgrade \par vrecica@matf.bg.ac.yu}

\end{document}